# Origin of the numerals

Ahmed Boucenna
Department of Physics, Faculty of Sciences,
Ferhat Abbas University 19000 Sétif, Algeria
aboucenna@wissal.dz   aboucenna@yahoo.com

**Abstract**
Through the pagination of an Arabian Algerian manuscript of the beginning of the 19th century, we rediscover the original shape, the "Ghubari" shape, of the numerals. Contrary to some assumptions, particularly those which claim that they are derived from Indian characters, this "Ghubari" shape, whose use has completely disappeared, shows that the ten modern numerals derive from ten Arabic letters. The symbol of a "Ghubari" numeral corresponds to the Arabic letter whose "Abjadi" numerical value is equal to this numeral. The assumption of the Indian origin of the numerals is denied by the shape of the numerals and by the right left sociological logic of the representation of the numerals and the algorithms of the basic operations.
The numerals are born in Maghreb or in Spain. In Europe, the "Ghubari" numerals became the modern numerals: 0, 1, 2, 3, 4, 5, 6, 7, 8, 9 and in the Middle East, borrowing two Hebrew letters, they gave the "Mashriki" numerals:   ٠   ١   ٢   ٣   ٤   ٥   ٦   ٧   ٨   ٩ .

MCS : 01A30

**Résumé**
A travers la pagination mixte d'un manuscrit arabe algérien du début du 19$^{ème}$ siècle nous redécouvrons la forme originale, la forme "Ghubari", des chiffres. Contrairement à ce que prétendent certaines hypothèses, particulièrement celles qui les présentent comme dérivant de caractères indiens, cette forme "Ghubari", dont l'utilisation a complètement disparu, montre que les dix chiffres modernes dérivent de dix lettres arabes. Le symbole d'un chiffre "Ghubari" correspond à la lettre arabe dont la valeur numérique "Abjadi" est égale à ce chiffre. L'hypothèse de l'origine indienne des chiffres et démentie par la forme des chiffres et par la logique sociologique droite gauche de la représentation des nombres et des algorithmes des opérations élémentaires.
Les chiffres sont nés au Maghreb ou en Espagne. En Europe, les chiffres "Ghubari" ont donné les chiffres modernes: 0, 1, 2, 3, 4, 5, 6, 7, 8, 9 et au Moyen Orient, en empruntant deux lettres hébraïques, ont donné les  chiffres "Mashriki" :   ٠   ١   ٢   ٣   ٤   ٥   ٦   ٧   ٨   ٩ .

## 1. Introduction
Before approaching the question of the origin of the modern numerals, let us first begin by recalling the following terminology. In this work, we designate by :
- **Modern numerals or Arabic modern numerals**, the numerals whose symbols are the following: 0, 1, 2, 3, 4, 5, 6, 7, 8, 9, which are used in Sciences and Technologies.
- **"Mashriki" numerals,** or **Arabic "Mashriki" numerals**, Whose symbols are :

٩ ٨ ٧ ٦ ٥ ٤ ٣ ٢ ١ ٠

which are currently used, up to now, in Middle East.
- **"Ghubari"** numerals, which are the ancestors of the modern numerals, were often used, during the Middle Ages, in Maghreb (North Africa) and Spain, at the time when the Arabic Muslim civilization was flourishing.

In this work, we approach the issue of the origin of the modern numeral symbols and that of the "Mashriki" ones. Some assumptions have been put forward as to this origin. One of these assumptions considers the modern numerals as deriving from Indian characters ([Ifrah, 1996], [Ifrah, 1999], [Ouaknin, 2004] and [Institut du monde arabe, 2005]). Another one makes a direct relationship between the number of plausible angles in the geometrical numeral shape and the numeral numerical value [Duchenoud, 1867]. The numeral would quantify the number of angles of its geometric shape. These assumptions seem to have as a starting point the modern shapes of the numeral symbols, whereas the origin must be searched through the original shapes that the numeral symbols had in Maghreb, before passing in Europe. They were then called the "Ghubari" numerals.

A first tentative to introduce the "Ghubari" numerals in Europe through Cordoba was undertaken by the Pope Sylvestre II (999-1003). In fact, it is from Béjaïa (Algeria) that the "Ghubari" numerals passed in Europe with tradesmen and Leonardo Fibonacci (born in 1170). In Europe, the evolution of the "Ghubari" numerals gave the modern ones, but in Maghreb they continued to be used under their original shape, regardless of their evolution in Europe until the beginning of the 19th century. Since the beginning of the second half of the 19th century, the use of the "Ghubari" numeral, in its original shape, completely disappeared, even in Maghreb, to the benefit of the "Mashriki" and modern numerals.

In this work, we first rediscover and identify the symbols of the "Ghubari" numerals through the mixed pagination of an Arabian Algerian manuscript, "Kitab khalil bni Ishak El Maliki", of the beginning of the 19th century. The correspondence, without ambiguity, between the modern numerals and their eldest the "Ghubari" numerals will be established. The elements of the "Abjadi" calculation ("Hissab El-joummel" or "Guematria"), particularly the notions of the numeric values ("Abjadi" numerical values) of the Arabic and Hebrew letters, will also be recalled. Then we will show the relation between the "Ghubari" numeral and the Arabic letter whose numerical value is equal to this numeral and we give the initial transformations imposed to the Arabic letters to become "Ghubari" numerals. We finish by unveiling the origin of the "Ghubari" numerals, the modern and the "Mashriki" ones. The right left sociological logic reinforced our findings.

## 2. Identification of the "Ghubari" numerals
### 2.1. Presentation of the manuscript "Kitab khalil bni Ishak El Maliki"
The manuscript is a copy of a book known as "Kitab khalil bni Ishak El Maliki", about the ritual Islamic Malikite "fik'h", prepared for "Sidi Khelifa" Zaouïa, in the region of Setif (Eastern Algeria.) The copyist is Mohamed ben Sidi khelifa ben Ali ben Salem ben Sifa. The manuscript is of 17.5 x cm 12.5 cm format and includes 179 sheets numbered from 1 to 179. These numbers were written, as a matter of fact, with the help of the "Ghubari" numerals. The sheets 4 and 23 are missing. The sheet number 77 has been numbered twice, having also number 78 on its back. The date of the copy, written with the "Maghribi" police of the Arabic letters, appears at the end of the manuscript, by means of the traditional expression" Copy finished at Midday ("Ezzaouel") of Wednesday, eight of the "Safar" month of the year 1225 of the "hegire" era, corresponding roughly to the year 1810, Figure 1. The year of the manuscript is indicated with "Ghubari" numerals above the expression عـــام, which means year. A transcription attempt, of this date in "Mashriki" numeral , which came later, can be seen as being unfruitful, since the numeral 5 has been taken for a 6 "Mashriki".

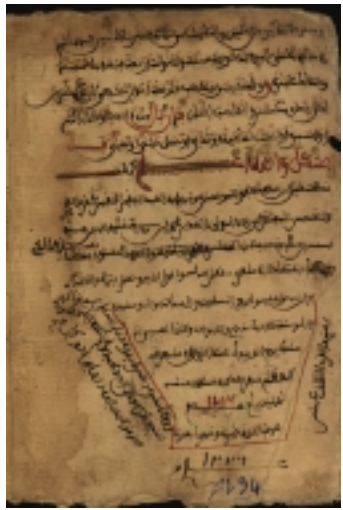

**Figure 1:**
Last page of the Manuscript
Kitab Khalil bni Ishak El Maliki

Another transcription was attempted, in modern numerals at the bottom of the page this time, also failed since the numeral 5 of the manuscript has been interpreted as being the numeral 4. This shows the total confusion to which people, who had just lost their "Ghubari" numerals for the benefit of the "Mashriki" and the modern ones, were confronted in the Maghreb during the second half of the 19th century.

### 2.2. The mixed pagination
Contrary to the pagination of modern books, it is rare to find, in the Arabian manuscripts, the order of the pages indicated by a series of numbers. It is the traditional pagination that was extensively used. The latter one consisted in writing down, in the left foot of even page of the sheet, the first word or the first two words, if there is ambiguity, of the first line of the odd following page, (see Figure 2).

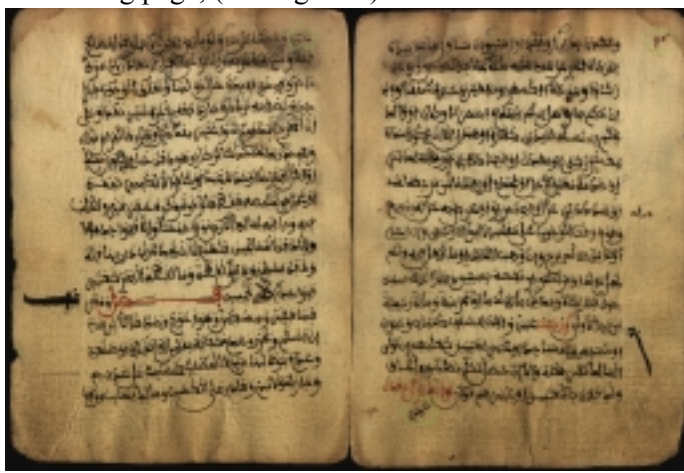

**Figure 2 :**
**Mixed pagination**
In top of the right page, is indicated the number of the sheet : 36.
At the bottom, on the left of the same page, the word عين is written which is the first word of the first line of the following page, the left page.

In our manuscript, in addition to the traditional pagination, the sheets have been numbered with the help of the "Ghubari" numerals. So, there is cohabitation between the traditional pagination and the one carried out with the help of the "Ghubari" numerals. It is the mixed pagination in which the order of the pages is reported in two different manners, on every sheet.

## 2.3. Identification of the symbols of the "Ghubari" numerals

The interest of the mixed pagination is double. First, it gives the symbols of the "Ghubari" numerals. Second, it eliminates all risk of mistake in the interpretation and the identification of the numeral symbols, since the number of numbered manuscript sheets is 179 and therefore, every "Ghubari" numeral is mentioned in more than 18 ordered numerical series. Figures 3 to 12 show the symbols of the "Ghubari" numerals used for the numbering of the sheets of the manuscript.

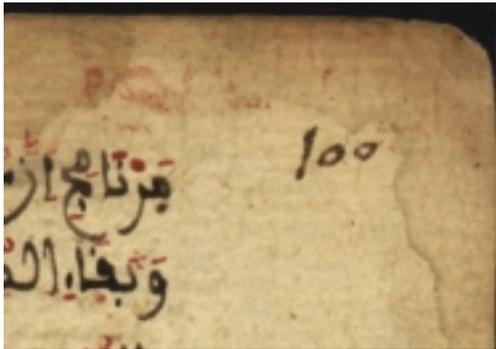

**Figure 3 :** from left to right
The "Ghubari" numerals 1, 0 and 0 in the manuscript.

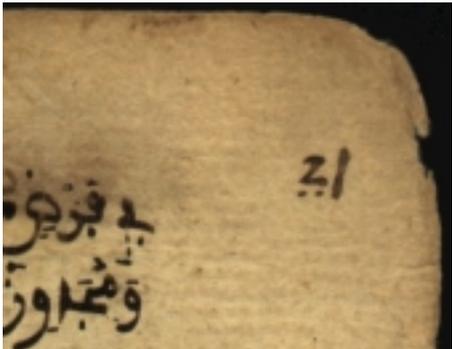

**Figure 4 :** from left to right
The "Ghubari" numerals 2 and 1 in the manuscript.

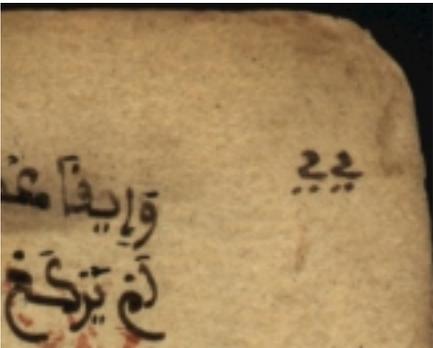

**Figure 5 :** from left to right
The "Ghubari" numerals 2 and 2 in the manuscript.

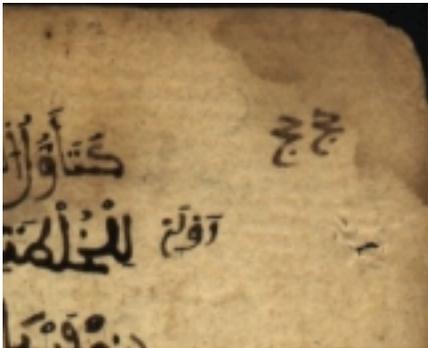

**Figure 6 :** from left to right
The "Ghubari" numerals 3 and 3 in the manuscript.

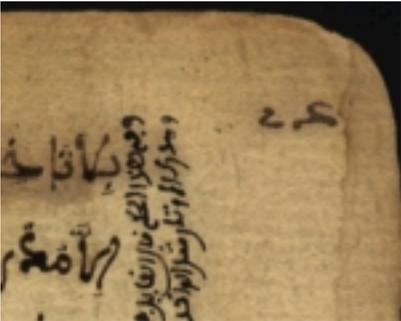

**Figure 7 :** from left to right
The "Ghubari" numerals 2 and 4 in the manuscript

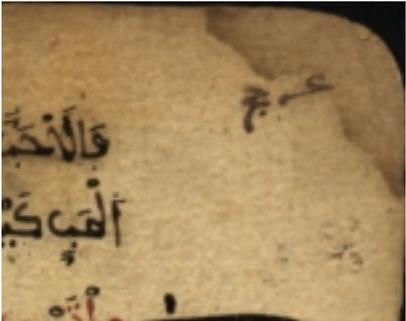

**Figure 8 :** from left to right
The "Ghubari" numerals 3 and 4 in the manuscript

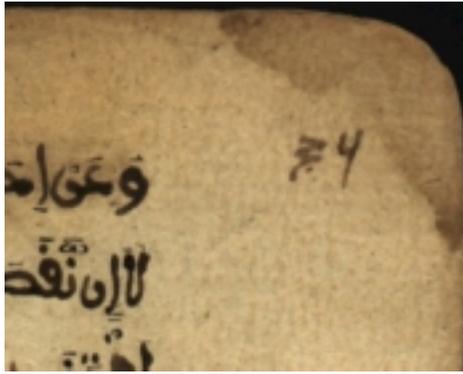

**Figure 9 :** from left to right
The "Ghubari" numerals 3 and 5 in the manuscript.

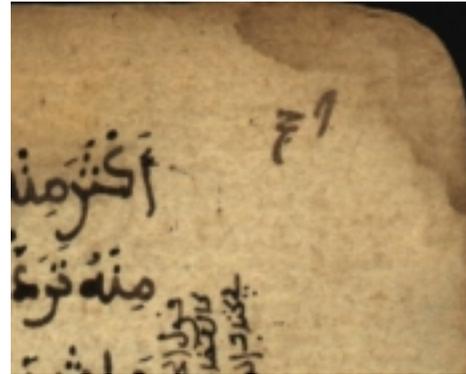

**Figure 10 :** from left to right
The "Ghubari" numerals 3 and 7 in the manuscript.

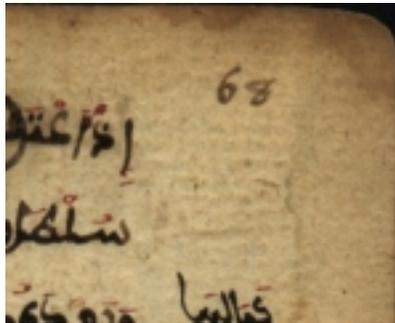

**Figure 11 :** from left to right
The "Ghubari" numerals 6 and 8 in the manuscript

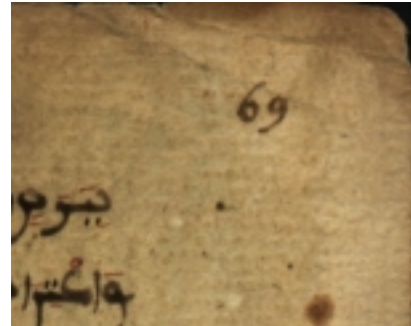

**Figure 12 :** from left to right
The "Ghubari" numerals 6 and 9 in the manuscript

The identification of "Ghubari" numeral symbols, deduced from the mixed pagination of the manuscript, is reported in Table 1. The correspondence is immediate between the modern numerals and their eldest the "Ghubari" numerals.

**Table 1**: The "Ghubari" numerals used in the numbering of the sheets
of the manuscript "Kitab khalil bni Ishak El Maliki"

| Modern numerals | 0 | 1 | 2 | 3 | 4 | 5 | 6 | 7 | 8 | 9 |
|---|---|---|---|---|---|---|---|---|---|---|
| "Ghubari" numerals | O | l | ݅ | ݇ | ے | 4 | 6 | ) | 8 | 9 |

The modern numerals symbols 0 and 1 are practically the same as the "Ghubari" symbols. In the modern numeral symbol 2, the two dots completely disappeared. In the modern numeral symbol 3, the dot and the tail disappeared. The comparison of the modern numerals 4 and 5 with their corresponding "Ghubari" ones shows that the two symbols have been permuted. Indeed, the numeral "Ghubari" 5 corresponds to the modern numeral 4, while the numeral "Ghubari" 4 corresponds to the modern numeral 5. The modern numeral symbol 6 is practically the same as the "Ghubari" numeral symbol 6. A slight evolution is observed in the modern numeral symbol 7. In the modern numeral symbol 8 the angles are rounded. The modern numeral symbol 9 does not differ from the numeral "Ghubari" symbol 9. The disappearance of the dots of the numerals 2 and 3 and the permutation of the numerals 4 and 5 constitute important divergences which, in our opinion, hide the origin of the numerals and mislead the researchers.

## 3. Elements of the "Abjadi" calculation
### 3.1. The "Abjadi" numerical values of the Arabic and Hebrew letters

The "Abjadi" calculation or "Hissab El-joummel" or Guématria, represents the set of the numerical operations that are possible from the numerical equivalence of the Arabic and Hebrew alphabetic letters. In the "Abjadi" calculation, the "Abjadi" numbers are the following natural numbers : {1, 2, 3, 4, 5, 6, 7, 8, 9, 10, 20, 30, 40, 50, 60, 70, 80, 90, 100, 200, 300, 400, 500, 600, 700, 800, 900, 1000}. Every Arabic and Hebrew letter possesses a numerical value ("Abjadi" numerical value) equal to an "Abjadi" number.

**Table 2**
Sounds and "Abjadi" numerical values of the 28 Arabic letters arranged according to the "Abjadi" order

| "Abjadi" Order | Arabic Letter | Sound | "Abjadi" numerical value | "Abjadi" Order | Arabic Letter | Sound | "Abjadi" numerical value |
|---|---|---|---|---|---|---|---|
| 1 | ا | Alif | 1 | 15 | س | Sin | 60 |
| 2 | ب | Baa | 2 | 16 | ع | A'in | 70 |
| 3 | ج | Jim | 3 | 17 | ف | Faa | 80 |
| 4 | د | Del | 4 | 18 | ص | Sad | 90 |
| 5 | ه | Haa | 5 | 19 | ق | K'af | 100 |
| 6 | و | Waw | 6 | 20 | ر | Raa | 200 |
| 7 | ز | Zin | 7 | 21 | ش | Shin | 300 |
| 8 | ح | H'aa | 8 | 22 | ت | Taa | 400 |
| 9 | ط | T'aa | 9 | 23 | ث | Thaa | 500 |
| 10 | ي | Yaa | 10 | 24 | خ | Kh'aa | 600 |
| 11 | ك | Kef | 20 | 25 | ذ | Dhel | 700 |
| 12 | ل | Lem | 30 | 26 | ض | Dzad | 800 |
| 13 | م | Mim | 40 | 27 | ظ | Dzaa | 900 |
| 14 | ن | Noun | 50 | 28 | غ | Ghin | 1000 |

**Table 3**
Sounds and "Abjadi" numerical values of the 22 Hebrew letters arranged according to the "Abjadi" order

| "Abjadi" Order | Hebrew Letter | Sound | "Abjadi" numerical value | "Abjadi" Order | Hebrew Letter | Sound | "Abjadi" numerical value |
|---|---|---|---|---|---|---|---|
| 1 | א | Aleph | 1 | 12 | ל | Lamedh | 30 |
| 2 | ב | Beth | 2 | 13 | מ ם | Mem | 40 |
| 3 | ג | Gimel | 3 | 14 | נ ן | Nun | 50 |
| 4 | ד | Daleth | 4 | 15 | ס | Samek | 60 |
| 5 | ה | He | 5 | 16 | ע | Ayin | 70 |
| 6 | ו | Vav | 6 | 17 | פ ף | Fe | 80 |
| 7 | ז | Zayin | 7 | 18 | צ ץ | Tsahde | 90 |
| 8 | ח | Cheth | 8 | 19 | ק | Q'oph | 100 |
| 9 | ט | Teth | 9 | 20 | ר | Regh | 200 |
| 10 | י | Yodh | 10 | 21 | שׂ שׁ | Sin Shin | 300 |
| 11 | כ ך | Kaph | 20 | 22 | ת | Tav | 400 |

The alphabet of a language is expressed according to a conventional order. In the "Abjadi" order, the Arabic and Hebrew letters are arranged in increasing order of their "Abjadi" numerical values. Table 2 gives the sounds and the "Abjadi" numerical values of the 28 Arabic letters arranged in the "Abjadi" order and Table 3 gives the sounds and the "Abjadi" numerical values of the 22 Hebrew letters arranged in the "Abjadi" order.

### 3.2. Representation of Numbers
Because of their "Abjadi" numerical values, the Arabic and Hebrew letters served as numerals to represent the numbers. All numbers between 1 and 1999 can be represented by a set of Arabic letters (a word), using the "Abjadi" numerical values of the letters. In the same way all numbers between 1 and 499 can be represented by a set of Hebrew letters (a word). To represent the number 1245, for example, the Arabic used the word "همرغ" composed of the Arabic letters: هـ of "Abjadi" numerical value 5, م of "Abjadi" numerical value 40, ر of "Abjadi" numerical value 200 and غ of "Abjadi" numerical value 1000. Notice the attached feature of the Arabic letters. Beyond, the class of the units is represented by a word having at most three Arabic letters, the class of the thousands is represented by another word, followed by the word "alf" (one thousand) to specify that it is thousands, the millions are represented by another word followed by the word "alf alf" to specify that it is millions (one thousand thousands). The number 23 456 789, for example, is written as:

<div align="center">طفذ و ونث ( ألف)  و  جك  ( ألف ألف)</div>

By taking into account the "Abjadi" numerical values of each letters, it reads as, [9 and 80 and 700] and [6 and 50 and 400] "alf" (thousand)  and [3 and 20] "alf alf" (thousand thousands : million)

### 3.3. Evaluation of expressions or "Hisseb el-joummel"
A whole sentence can be evaluated while calculating the sum of the "Abjadi" numerical values of the letters composing the words of this sentence. For example, the value of the expression "احمد زينب" is:
$$1 + 8 + 40 + 4 + 7 + 10 + 50 + 2 = 122.$$
This kind of digitalisation of expressions using the "Abjadi" numerical values of the letters is used in astrology. The values of the expressions are treated with well defined algorithm calculations and the final result is compared to astrological tables.

## 4. "Ghubari" numerals and Arabic letters
### 4.1. The "Ghubari" numerals 1, 2, 3, 4, 5, 6, 7, 8 and 9
The dot and the shape of the "Ghubari" numeral 3, do recall us the Arabic letter ج (Jim), and one is attempted to ask the following question:

*Are the "Ghubari" numerals simply Arabic letters?*

Already, the fact that the "Ghubari" numeral 4 is our modern numeral 5 proves that the hypothesis that makes a direct relation between the number of possible angles in the geometrical numeral shape and the numeral value is false, since the "Ghubari" numerals 4 and 5 do not correspond to the number of angles of their geometrical shapes.
In Table 4, we show the first ten Arabic letters of the "Abjadi" order written with the "Maghribi" police and the "Ghubari" numerals as identified from our manuscript. Let us now examine in more details, one by one, all the symbols of the "Ghubari" numerals and let us compare them, to the symbols of the first ten Arabic letters of the "Abjadi" order.
The Arabic letter ج (Jim) has indeed the "Abjadi" numerical value 3, correctly corresponding to the numeral "Ghubari" 3, as we suspected it. The coincidence is remarkable. For all "Ghubari numerals, we can note that, (see Table 4):
The "Ghubari" numeral 1 precisely corresponds to the Arabic letter ا (Alif), of "Abjadi" numerical value 1.

The "Ghubari" numeral 2 with its two dots, does not have any tie with the Arabic letter ب (Baa) of "Abjadi" numerical value 2. The two dots lead us to think about the Arabic letter ي (Yaa) of "Abjadi" numerical value 10. It is precisely one of the final shape of the Arabic letter ي (Yaa) written in the "Maghribi" police. In the manuscript, the two dots are sometimes omitted as with the two dots of the final Arabic letter ى (Yaa).

The "Ghubari" numeral 3 is the slightly modified final bounded shape of the Arabic letter ج (Jim) of "Abjadi" numerical value 3. The presence of the dot of the ج (Jim) in the "Ghubari" numeral 3 can be noticed. The letter ج (Jim) has been doubled by its initial shape, without the dot however.

**Table 4**
The "Ghubari" numerals, the first ten Arabic letters of the "Abjadi" order
and the transformations imposed to the Arabic letters to become "Ghubari" numerals

| Arabic letter | Sound | "Abjadi" numerical value | "Maghribi" shaped arabic leter | Transformation | "Ghubari" numeral |
|---|---|---|---|---|---|
| ا | Alif | 1 | ا | NONE | 1 |
| ب | Baa | 2 | ب | | |
| ي | Yaa | 10 | ي ي | NONE | 2 |
| ج | Jim | 3 | ج ج | ج | 3 |
| د | Del | 4 | د | NONE | 4 |
| ه | Haa | 5 | ه ه ه ه | 4 | 4 |
| و | Waw | 6 | و | UP-SIDE-DOWN | 6 |
| ز | Zin | 7 | ز ن | ر | 7 |
| ح | H'aa | 8 | ح | 8 | 8 |
| ط | T'aa | 9 | ط | UP-SIDE-DOWN | 9 |
| ص | Sad | 90 | ص ص | NONE | 0 |

The "Ghubari" numeral 4 is precisely the "Maghribi" shape of the Arabic letter د (Del), of "Abjadi" numerical value 4.

The "Ghubari" numeral 5 is the modified final bounded shape of the Arabic letter ه (Haa), of "Abjadi" numerical value 5.

The "Ghubari" numeral 6 is the up-side down shape of the Arabic letter و (Waw), of "Abjadi" numerical value 6.

The "Ghubari" numeral 7 is the slightly modified shape of the Arabic letter ز (Zin), of "Abjadi" numerical value 7. The dot has been substituted by a small feature linked to the body of the letter.

The "Ghubari" numeral 8 is the final shape of the Arabic letter ح (H'aa), of "Abjadi" numerical value 8. The tail of the letter ح (H'aa) is bound to its beginning.

The "Ghubari" numeral 9 is the up-side down shape of the Arabic letter ط (Taa), of "Abjadi" numerical value 9.

The correspondence between the shape of the Arabic letter whose "Abjadi" numerical value is equal to the "Ghubari" numeral and the symbol of this numeral is perfect for the "Ghubari" numerals 1, 3, 4, 5, 6, 7, 8 and 9. The coincidences are remarkable. The "Ghubari" numeral 2 does not correspond to the second Arabic letter of the "Abjadi" order but corresponds to the Arabic letter having the "Abjadi" numerical value 10.

The "Ghubari" numerals, ancestors of our modern numerals, are therefore derived from the first ten Arabic letters of the "Abjadi" order. The initial strategy for obtaining the "Ghubari" numeral symbols consisted therefore to choose, as a symbol of a "Ghubari" numeral, the Arabic letter whose "Abjadi" numerical value is equal to this numeral.

**4.2. The "Ghubari" numeral 0**

The discovery of 0, in its present conception which has simplified the representation of the numbers and the algorithms of the basic operations, was certainly a very crucial event. If the choice of the symbols representing the "Ghubari" numerals 1, 2,…, 9 can be justified by the "Abjadi" numerical values of the Arabic letters chosen to represent them, the symbol representing the "Ghubari" numeral 0 cannot be explained in the same way since the value zero is not an "Abjadi" numerical value. Like all new inventions, it is necessary to choose for it a name and a symbol.

The chosen name is the word "Sifr". Is it an Arabic word? Figure 13 shows a part of an Arabic dictionary page which gives the meaning of the word "Sifr" and its derivatives. One can see that the word "Sifr" is not a foreign word to Arabic. It has not been borrowed from another language to describe a new state, unknown in Arabic, as for example the case of the word "Falsafa" that refers to "philosophy", which is borrowed from Greek to describe a new state. The word "Sifr" does not derive, as some authors,([Ifrah, 1996], [Ifrah, 1999], [Ouaknin, 2004] and [Institut du monde arabe, 2005]) pretend, from of the Indian word "Shûnya", since the word "Sifr" and its derivatives existed in Arabic long before the appearance of zero itself.

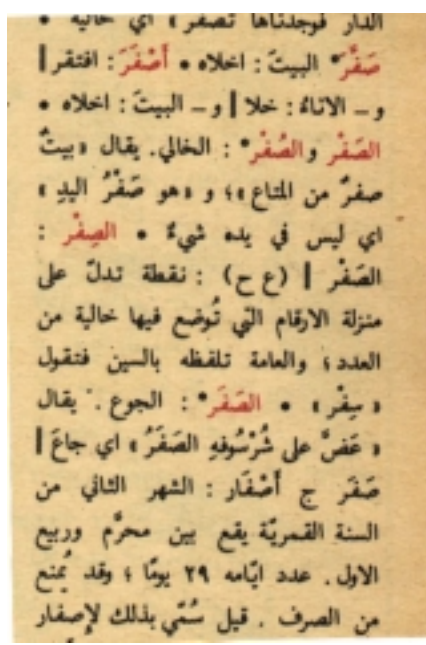

**Figure 13** :
Meaning of the Arabic word "Sifr"
and of its derivatives

In Arabic the meaning of the word "sifr" describes an empty state which one did not expect. The emptiness described by the word "sifr" is in all cases an abnormal situation. The name "Sifr" given to 0 expresses exactly the role that the 0 must play in the representation of the numbers with the help of the numerals. The numeral 0 must replace a numeral that made defect in a given rank. It does not represent the emptiness, it fills the emptiness.

The chosen symbol to represent the "Sifr" is the symbol 0. Does this symbol derive from an Arabic letter as the other "Ghubari" numerals? If yes, which one? As the value 0 is not an "Abjadi" numerical value, one cannot find an Arabic letter having the good "Abjadi" numerical value, to represent the "Ghubari" numeral 0. The strategy of choice of the symbols must be modified so that one will not take into account the "Abjadi" numerical values of the letters. One could think then about choosing the tenth Arabic letter ي (Yaa) of "Abjadi" numerical value 10 to represent the "Ghubari" numeral 0. This was not the case. In my opinion, it is the initial shape of the Arabic letter ص (Sad), first letter of the Arabic word صفر "Sifr", which served to manufacture the symbol of the "Ghubari" numeral 0. The choice of the symbol 0, has not been motivated therefore as some authors pretend ([Ifrah, 1996], [Ifrah, 1999], [Ouaknin, 2004] and [Institut du monde arabe, 2005]), by high philosophical considerations. The fact that the difference of two equal numbers is exactly equal to the numeral 0, is a consequence of the representation of the numbers with the help of the "Ghubari" numerals.

**4.3. The initial transformations**

The initial strategy for obtaining the "Ghubari" numeral symbols consisted to choose, as a symbol of a "Ghubari" numeral, the Arabic letter whose "Abjadi" numerical value is equal to this numeral. As the numerical value 0 is not an "Abjadi" numerical value, the numeral "Ghubari" 0 was the first numeral formed regardless of the "Abjadi" numerical values of the Arabic letters. The idea to manufacture the other "Ghubari" numerals 1, 2,…, 9 independently of the Arabic letters comes right after to confirm definitively the separation between the numerals and the letters with their "Abjadi" numerical values. Indeed, the "Abjadi" numerical values would have led to the representation of the "Ghubari" numeral 5 by the isolated shape of the letter ه (Haa), having the "Abjadi" numerical value 5. However this kind of circle has probably already been used to designate the numeral "Ghubari" 0. It was necessary to find another symbol to designate the numeral "Ghubari" 5. One deduces from this priority that the "Ghubari" numeral 0 was the first "Ghubari" numeral put in shape. The table 5 gives the likely initial transformations having succeeded to the "Ghubari" numerals. In the likely initial version of the "Ghubari" numerals, the numerals are simply Arabic letters, the "Abjadi" order of the Arabic letters was scrupulously respecting. In the Final version of the "Ghubari" numerals, the initial transformations imposed by the Master of the numeral were taking into account. The invention of the "Ghubari" numerals has been achieved in a very short time. Once the idea of the numeral discovered, it took 1 to 2 days to form all "Ghubari" numerals including the 0. This explains the freedom taken by the Master of the numerals, during the elaboration of all the "Ghubari" numerals. He did not take the Arabic letter ب (Baa), according to the "Abjadi" order of the letters to designate the numeral "Ghubari" 2, but he chose the final shape of the "Maghribi" shaped Arabic letter ي (Yaa) with its two dots. Arabic words are written with attached letters, to put in shape the "Ghubari" numerals, the Master made some transformations in the Arabic letters in order to hide their attached feature. So he did not take directly the letter و (Waw) of "Abjadi" numerical value 6 to represent the "Ghubari" numeral 6, but he decided to take the up-side down shape of the letter و (Waw). He did not take directly the letter ط (Taa) of "Abjadi" numerical value 9 to represent the numeral 9, but he took the up-side down shape of the letter ط (Taa). Our Master had the same attitude during the formation of all "Ghubari" numerals.

**Table 5**
Likely initial transformations having succeeded to the "Ghubari" numerals and the modern numerals

| Arabic letter | Sound | abjadi numerical value | Initial version of the "Ghubari" Numeral | final version of the "Ghubari" Numeral | European version of "Ghubari" Numeral | Chiffres Modernes |
|---|---|---|---|---|---|---|
| ا | Alif | 1 | ا | ا | ا | 1 |
| ب | Baa | 2 | ب | ح | ح | 2 |
| ج | Jim | 3 | ج | ج | ج | 3 |
| د | Del | 4 | د | ح | 4 | 4 |
| ه | Haa | 5 | ه | 4 | ح | 5 |
| و | Waw | 6 | و | 6 | 6 | 6 |
| ز | Zin | 7 | ز | 7 | 7 | 7 |
| ح | H'aa | 8 | ح | 8 | 8 | 8 |
| ط | T'aa | 9 | ط | 9 | 9 | 9 |
| ي | Yaa | 10 | | | | |
| ص | Sad | | ص | O | O | 0 |

This freedom shows the exhibited intention to completely depart from the old representation of the numbers that was intimately bounded to the "Abjadi" numerical values of the Arabic letters. It is a new era that begins and our Master was well conscious of the importance of his realization. He acted freely, with a new spirit, but without completely forgetting the "Abjadi" numerical values of the Arabic letters. Therefore, the shapes of the "Ghubari" numerals really originated from the Arabic letters but their spirit was completely different from the spirit of the "Abjadi" numerical values.

**4.4. The right left sociological Logic**
The numerals have an interest by their impacts on the representation of the numbers and on the algorithm calculations. When closely examining this representation of the numbers and these algorithm calculations, we discover a logic, called here, the right left sociological logic.
If a person has the right left sociological logic, Figure 14, he writes and reads a word or a number from the right to the left. To write, for example, the number : 12457892, he starts by writing the units, the ten, the hundreds of the units, then the units, the ten, the hundreds of the thousands etc. starting by the right. The same person will proceed similarly for reading: starting by right: 2 and 90 and 800; 7 and 50 and 400 thousands ; 2 and 10 million. In the left right sociological logic, one writes and reads a word or a number from left to right. Reading the same number 12457892 which is written with the right left sociological logic, one starts by putting himself in the right-left sociological logic by isolating the numerals composing this number from right to left as follow : 12 457 892 and returns to his left right sociological logic, and reads : 12 million 457 thousands 892.

**Figure 14**: Sociological Logic

Reading the number: 12457892 ?

A person who has the Right Left Sociological Logic :

Reading in one step :

2 and 90 and 800 and 7 and 50 and 400 thousand and 2 and 10 millions

A person who has the Left Right Sociological Logic:

Reading in two steps :

First step : using the Right Left Sociological Logic:

12 457 892

Second Step : returning to the Left Right Sociological Logic and reading:

12 millions 457 one thousand 892

All basic operations, for which "Ghubari" numerals simplified the algorithms, were accomplished in the right left sociological logic. This logic is imposed to us, by the "Ghubari" numerals, every time we deal with numbers. When stocking a number, in the computer memory, for example, the weakest bit is in the right position.

**5. Evolution of the numeral**
The "Ghubari" numerals are born in Maghreb or in Spain, emigrated toward the Mashrik to become, after evolution, the "Mashriki" numerals and toward Europe to become after evolution, the modern numerals.

**5.1. The modern numerals**
Except the permutation of the numerals 4 and 5 and the disappearance of the dots of the numerals 2 and 3, the comparison of the "Ghubari" numerals and modern numerals shows that the transformations especially touched the esthetic form. Table 5 shows the evolution of the "Ghubari" numerals to European version of "Ghubari" numerals and toward the modern numerals.

**5.2. The "Mashriki" numerals**
In the Middle East, the "Mashriki" numerals are used up to now. From where do these "Mashriki" numerals come ? The meticulous observation of the "Mashriki" numerals, particularly the numerals 1, 4, and 9 permits to affirm that they have the "Ghubari" numerals, and therefore the Arabic letters, as origin, see Table 6. The changes made to pass for "Ghubari" numerals to the "Mashriki" numerals likely were meant to re-establish the ascendancy of the "Abjadi" numerical values of the Arabic letters on the numerals. The strategy of the transformations consisted to:
- pratically redraw the shape of the numerals by scrupulously respecting the "Abjadi" order of the letters, starting with the numeral 1. The numeral 0 is not the first numeral put in shape,
- replace the dots of the Arabic letters by a sort of leg attached to the body of the letter,
- avoid taking identical symbols to the "Ghubari" numerals to represent different numerals to prevent all confusion,
- borrow Hebrew letters instead of Arabic ones every time that might be confusion with the "Ghubari" numerals.

The table 6 shows the transformations undergone on the Arabic and Hebrew letters and the "Ghubari" numerals to get "Mashriki" numerals. The likely initial version of the "Mashriki" numerals respects scrupulously the Arabic letters "Abjadi" order and the final version of the "Mashriki" numerals, takes into account the mentioned transformations.
One notes that:
The "Mashriki" numeral 1: is precisely the Arabic letter ا (Alif), of the "Abjadi" numerical value 1.
The "Mashriki" numeral 2: is represented by a modified shape of the Arabic letter ب (Baa), of "Abjadi" numerical value 2. The dot is replaced by a sort of a leg attached to the body of the letter. The numeral 2 is not represented by the Arabic letter ي (Yaa) corresponding to the "Abjadi" numerical value 10, as was the case for the "Ghubari" numeral 2.
The "Mashriki" numeral 3: is a shape; modified by a rotation of π/2 to the right, of the numeral "Ghubari" 3, wich is the slightly modified final bound shape of the Arabic letter ج (Jim) of "Abjadi" numerical value 3. As for the numeral 2, the dot is replaced by a leg attached to the body of the letter.
The "Mashriki" numeral 4: is a slightly distorted shape of the "Ghubari" numeral 4, which is precisely the "Maghribi" shape of the Arabic letter د (Del), of "Abjadi" numerical value 4.

**Table 6**
Evolution of the Arabic and Hebrew letters and "Ghubari" numerals toward the "Mashriki" numerals

| Abjadi Numerical value | Arabic Letter | Sound | Hebrew letter | Sound | Transformation | Likely Initial version of the Mashriki numeral | Final version of the Mashriki numeral | Ghubari numeral |
|---|---|---|---|---|---|---|---|---|
| 1 | ا | Alif | א | Aleph | NONE | ا | ١ | ١ |
| 2 | ب | Baa | ב | Beth | ݑ | ب | ٢ | |
| 3 | ج | Jim | ג | Gimel | | ج | ٣ | |
| 4 | ح | Del | ד | Daleth | NONE | ح | ٤ | |
| 5 | ه | Haa | ה | He | NONE | ه | ٥ | 4 |
| 6 | و | Waw | ו | Vav | NONE | و | ٦ | 6 |
| 7 | ز | Zin | ז | Zayin | ں | ز | ٧ | ɿ |
| 8 | ح | H'aa | ח | Cheth | ROTATION π/2 | ح | ٨ | 8 |
| 9 | ط | T'aa | ט | Teth | UP-SIDE-DOWN | ط | ٩ | 9 |
| 10 | ي | Yaa | י | Yodh | NONE | ي | ٠ | O |

The "Mashriki" numeral 5: is represented by the isolated shape of the Arabic letter ه (Haa), of "Abjadi" numerical value 5. This shape is very close to that of the one kept by the Master of the numeral to represent the "Ghubari" numeral 0.

The "Mashriki" numeral 6: has been represented by the shape of the Hebrew letter ו (Vav), of "Abjadi" numerical value 6. The choice of the Arabic letter و (Waw) of "Abjadi" numerical value 6 is prevented by the possible confusion with the "Ghubari" numeral 9.

The "Mashriki" numeral 7: is a modified shape of the Arabic letter ز (Zin), of "Abjadi" numerical value 7. Like the case of the numeral 2, the dot is replaced by a leg attached to the body of the letter.

The "Mashriki" numeral 8 : is a modified shape, by a rotation of π/2 to the left, of the Arabic letter ح (H'aa), of "Abjadi" numerical value 8.

The "Mashriki" numeral 9 : is the up-side down shape of the Arabic letter ط (Taa), of "Abjadi" numerical value 9. It is also the symbol of the "Ghubari" numeral 9.

The "Mashriki" numeral 8 : is a modified shape, by a rotation of π/2 to the left, of the Arabic letter ح (H'aa), of "Abjadi" numerical value 8. The shape of The "Mashriki" numeral 8 lead us to think about the hebrew letter ח(cheteh) of "Abjadi" value 8

The "Mashriki" numeral 9 : is the up-side down shape of the Arabic letter ط (Taa), of "Abjadi" numerical value 9. It is also the symbol of the "Ghubari" numeral 9.

The "Mashriki" numeral 0 : has been represented by the shape of the Hebrew letter the letter י (Yodh), of "Abjadi" numerical value 10. The symbol 0 kept by the Master of the numeral to represent the "Ghubari" numeral 0 has a shape very close to the one ه (haa) chosen to represent the numeral "Mashriki" 5. Therefore it was necessary to choose another symbol for the Mashriki" numeral 0. A possible strategy of choice of the symbols would consist in taking the Arabic letter ي (Yaa) of "Abjadi" numerical value 10 as symbol of the "Mashriki" numeral 0. This choice may lead to confusion with the numeral "Ghubari" 2. To prevent this confusion, it might be possible to borrow the Hebrew letter י (Yodh), of "Abjadi" numerical value 10. This choice shows that the first "Mashriki" numeral put in shape was not the numeral 0. The choice of the symbol "Mashriki" numeral 0, has not been motivated by high philosophical considerations.

## 6. CONCLUSION

Since his early beginning, the man tried to quantify the nature by using numbers. To represent these numbers he used numerals. Every civilization had its numerals. The "Ghubari" numerals, with their zero, were, for different reasons, a revolution. They simplified the representation of the numbers and the algorithm calculations. They provided us with a small calculator to perform complicated operations. This mini calculator is not an instrument but a simple method that only requires a sheet and a pencil.

Through the mixed pagination of an Arabian Algerian manuscript of the beginning of the 19th century we rediscover the "Ghubari" numerals whose use has completely disappeared. The "Ghubari" numerals gave in Europe the modern numerals. Contrary to some assumptions, the "Ghubari" shape shows that the ten numerals we presently use are ten slightly modified Arabic letters given in the "Abjadi" order. The symbol of a "Ghubari" numeral corresponds to the Arabic letter whose "Abjadi" numerical value is equal to this numeral.

The right left sociological logic is imposed to us, in the representation (reading and writing) of the numbers and in the algorithms of the basic operations, by the "Ghubari" numerals, every time we deal with numbers.

The "Ghubari" numerals have for origin the "Maghribi" shaped Arabic letters. The "Ghubari" numerals that simplified the writing of the numbers and the algorithms of the basic operations have been often used in the Ghubari (calculation) in the Maghreb until the beginning of the 19th century.

From what we have just seen we can conclude that, the ten "Ghubari" numerals are originated either in the Maghreb or in Spain during the flourishing time of the Arabic Muslim civilization.

From the city port Béjaïa (Algeria) the "Ghubari" numerals travelled and introduced to Europe and became there, after evolution the modern numerals: 0, 1, 2, 3, 4, 5, 6, 7, 8, 9. They also travelled to the Middle East to became, after transformations and by adding two Hebrew letters, the "Mashrikis" numerals: ٠ ١ ٢ ٣ ٤ ٥ ٦ ٧ ٨ ٩ .

For historical reasons, decline and colonization, the "Ghubari" numerals became, during the second half of the 19th century and the first half of the 20th century, foreign and unknown in their own country, the Maghreb. Today the "Ghubari" numerals, becoming the modern numerals, have returned from Europe to be the official numerals in the Maghreb countries (Algeria, Morocco and Tunisia).